\documentstyle[twoside]{article}

\catcode`\@=11
\long\def\@makefntext#1{
\protect\noindent \hbox to 3.2pt {\hskip-.9pt
$^{{\eightrm\@thefnmark}}$\hfil}#1\hfill}               %CAN BE USED

\def\thefootnote{\fnsymbol{footnote}}
\def\@makefnmark{\hbox to 0pt{$^{\@thefnmark}$\hss}}    %ORIGINAL

\def\ps@myheadings{\let\@mkboth\@gobbletwo
\def\@oddhead{\hbox{}
\rightmark\hfil\eightrm\thepage}
\def\@oddfoot{}\def\@evenhead{\eightrm\thepage\hfil
\leftmark\hbox{}}\def\@evenfoot{}
\def\sectionmark##1{}\def\subsectionmark##1{}}

\oddsidemargin=\evensidemargin
\renewcommand{\thefootnote}{\fnsymbol{footnote}}

\newcounter{sectionc}\newcounter{subsectionc}\newcounter{subsubsectionc}
\renewcommand{\section}[1] {\vspace{12pt}\addtocounter{sectionc}{1}
\setcounter{subsectionc}{0}\setcounter{subsubsectionc}{0}\noindent
        {\tenbf\thesectionc. #1}\par\vspace{5pt}}
\renewcommand{\subsection}[1] {\vspace{12pt}\addtocounter{subsectionc}{1}
        \setcounter{subsubsectionc}{0}\noindent
        {\bf\thesectionc.\thesubsectionc. {\kern1pt \bfit #1}}\par\vspace{5pt}}
\renewcommand{\subsubsection}[1] {\vspace{12pt}\addtocounter{subsubsectionc}{1}
        \noindent{\tenrm\thesectionc.\thesubsectionc.\thesubsubsectionc.
        {\kern1pt \tenit #1}}\par\vspace{5pt}}
\newcommand{\nonumsection}[1] {\vspace{12pt}\noindent{\tenbf #1}
        \par\vspace{5pt}}

\newcounter{appendixc}
\newcounter{subappendixc}[appendixc]
\newcounter{subsubappendixc}[subappendixc]
\renewcommand{\thesubappendixc}{\Alph{appendixc}.\arabic{subappendixc}}
\renewcommand{\thesubsubappendixc}
        {\Alph{appendixc}.\arabic{subappendixc}.\arabic{subsubappendixc}}

\renewcommand{\appendix}[1] {\vspace{12pt}
        \refstepcounter{appendixc}
        \setcounter{figure}{0}
        \setcounter{table}{0}
        \setcounter{lemma}{0}
        \setcounter{theorem}{0}
        \setcounter{corollary}{0}
        \setcounter{definition}{0}
        \setcounter{equation}{0}
        \renewcommand{\thefigure}{\Alph{appendixc}.\arabic{figure}}
        \renewcommand{\thetable}{\Alph{appendixc}.\arabic{table}}
        \renewcommand{\theappendixc}{\Alph{appendixc}}
        \renewcommand{\thelemma}{\Alph{appendixc}.\arabic{lemma}}
        \renewcommand{\thetheorem}{\Alph{appendixc}.\arabic{theorem}}
        \renewcommand{\thedefinition}{\Alph{appendixc}.\arabic{definition}}
        \renewcommand{\thecorollary}{\Alph{appendixc}.\arabic{corollary}}
        \renewcommand{\theequation}{\Alph{appendixc}.\arabic{equation}}
        \noindent{\tenbf Appendix \theappendixc #1}\par\vspace{5pt}}
\newcommand{\subappendix}[1] {\vspace{12pt}
        \refstepcounter{subappendixc}
        \noindent{\bf Appendix \thesubappendixc. {\kern1pt \bfit #1}}
        \par\vspace{5pt}}
\newcommand{\subsubappendix}[1] {\vspace{12pt}
        \refstepcounter{subsubappendixc}
        \noindent{\rm Appendix \thesubsubappendixc. {\kern1pt \tenit #1}}
        \par\vspace{5pt}}

\newcommand{\textlineskip}{\baselineskip=13pt}
\newcommand{\smalllineskip}{\baselineskip=10pt}

\def\eightcirc{
\begin{picture}(0,0)
\put(4.4,1.8){\circle{6.5}}
\end{picture}}
\def\eightcopyright{\eightcirc\kern2.7pt\hbox{\eightrm c}}

\def\abstracts#1#2#3{{
        \centering{\begin{minipage}{4.5in}\baselineskip=10pt\footnotesize
        \parindent=0pt #1\par
        \parindent=15pt #2\par
        \parindent=15pt #3
        \end{minipage}}\par}}

\renewenvironment{thebibliography}[1]                   %ALL CHANGES DD 13/3/92
        {\frenchspacing
         \ninerm\baselineskip=11pt
         \begin{list}{\arabic{enumi}.}
        {\usecounter{enumi}\setlength{\parsep}{0pt}
         \setlength{\leftmargin 12.7pt}{\rightmargin 0pt} %FOR 1--9 ITEMS
         \setlength{\itemsep}{0pt} \settowidth
        {\labelwidth}{#1.}\sloppy}}{\end{list}}

\newcounter{itemlistc}
\newcounter{romanlistc}
\newcounter{alphlistc}
\newcounter{arabiclistc}

\newcommand{\fcaption}[1]{
        \refstepcounter{figure}
        \setbox\@tempboxa = \hbox{\footnotesize Fig.~\thefigure. #1}
        \ifdim \wd\@tempboxa > 5in
           {\begin{center}
        \parbox{5in}{\footnotesize\smalllineskip Fig.~\thefigure. #1}
            \end{center}}
        \else
             {\begin{center}
             {\footnotesize Fig.~\thefigure. #1}
              \end{center}}
        \fi}

\newcommand{\tcaption}[1]{
        \refstepcounter{table}
        \setbox\@tempboxa = \hbox{\footnotesize Table~\thetable. #1}
        \ifdim \wd\@tempboxa > 5in
           {\begin{center}
        \parbox{5in}{\footnotesize\smalllineskip Table~\thetable. #1}
            \end{center}}
        \else
             {\begin{center}
             {\footnotesize Table~\thetable. #1}
              \end{center}}
        \fi}

\def\@citex[#1]#2{\if@filesw\immediate\write\@auxout
        {\string\citation{#2}}\fi
\def\@citea{}\@cite{\@for\@citeb:=#2\do
        {\@citea\def\@citea{,}\@ifundefined
        {b@\@citeb}{{\bf ?}\@warning
        {Citation `\@citeb' on page \thepage \space undefined}}
        {\csname b@\@citeb\endcsname}}}{#1}}

\newif\if@cghi
\def\cite{\@cghitrue\@ifnextchar [{\@tempswatrue
        \@citex}{\@tempswafalse\@citex[]}}
\def\citelow{\@cghifalse\@ifnextchar [{\@tempswatrue
        \@citex}{\@tempswafalse\@citex[]}}
\def\@cite#1#2{{$\null^{#1}$\if@tempswa\typeout
        {IJCGA warning: optional citation argument
        ignored: `#2'} \fi}}

\def\pmb#1{\setbox0=\hbox{#1}
        \kern-.025em\copy0\kern-\wd0
        \kern.05em\copy0\kern-\wd0
        \kern-.025em\raise.0433em\box0}

\def\fnt#1#2{\footnotetext{\kern-.3em
        {$^{\mbox{\scriptsize #1}}$}{#2}}}

\def\fpage#1{\begingroup
\voffset=.3in
\thispagestyle{empty}\begin{table}[b]\centerline{\footnotesize #1}
        \end{table}\endgroup}

\def\runninghead#1#2{\pagestyle{myheadings}
\markboth{{\protect\footnotesize\it{\quad #1}}\hfill}
{\hfill{\protect\footnotesize\it{#2\quad}}}}
\font\tenrm=cmr10
\font\tenit=cmti10
\font\tenbf=cmbx10
\font\bfit=cmbxti10 at 10pt
\font\ninerm=cmr9
\font\nineit=cmti9
\font\ninebf=cmbx9
\font\eightrm=cmr8

\newtheorem{lemma}{Lemma}
\newtheorem{corollary}{Corollary}
\def\qed{\hbox{${\vcenter{\vbox{                        %HOLLOW SQUARE
   \hrule height 0.4pt\hbox{\vrule width 0.4pt height 6pt
   \kern5pt\vrule width 0.4pt}\hrule height 0.4pt}}}$}}

\renewcommand{\thefootnote}{\fnsymbol{footnote}}        %USE SYMBOLIC FOOTNOTE

\def\bsc{{\sc a\kern-6.4pt\sc a\kern-6.4pt\sc a}}       %LATEX LOGO
\def\bflatex{\bf L\kern-.30em\raise.3ex\hbox{\bsc}\kern-.14em
T\kern-.1667em\lower.7ex\hbox{E}\kern-.125em X}

\def\Mult#1#2{\left[{#1 \atop #2}\right]}

\begin{document}

\runninghead{Higher-level Bailey lemma}
 {Higher-level Bailey lemma}

\normalsize\textlineskip
\thispagestyle{empty}
\setcounter{page}{1}

%\copyrightheading{}                     %{Vol. 0, No. 0 (1993) 000---000}

\vspace*{0.88truein}

\fpage{1}
\centerline{\bf A HIGHER-LEVEL BAILEY LEMMA}
\vspace*{0.37truein}
\centerline{\footnotesize ANNE SCHILLING\footnote{
E-mail: anne@insti.physics.sunysb.edu}}
\vspace*{0.015truein}
\centerline{\footnotesize\it 
Institute for Theoretical Physics, State University of New York}
\baselineskip=10pt
\centerline{\footnotesize\it Stony Brook, NY 11794-3840, USA}
\vspace*{10pt}
\centerline{\normalsize and}
\vspace*{10pt}
\centerline{\footnotesize S.~OLE WARNAAR\footnote{
E-mail: warnaar@maths.mu.oz.au}}
\vspace*{0.015truein}
\centerline{\footnotesize\it 
Mathematics Department, University of Melbourne}
\baselineskip=10pt
\centerline{\footnotesize\it 
Parkville, Victoria 3052, Australia}
\vspace*{0.225truein}
\begin{center}\footnotesize\smalllineskip April 1996\\ \end{center}
\vspace*{0.21truein}
\abstracts{We propose a generalization of Bailey's lemma, 
useful for proving $q$-series identities.
As an application, generalizations of
Euler's identity, the Rogers--Ramanujan identities,
and the Andrews--Gordon identities are derived.
This generalized Bailey lemma also allows one to derive identities for the
branching functions of higher-level A$^{(1)}_1$ cosets.
}{}{}

\setcounter{footnote}{0}
\renewcommand{\thefootnote}{\alph{footnote}}

\vspace*{1pt}\textlineskip      %) USE THIS MEASUREMENT WHEN THERE IS
\section{The Bailey Lemma}      %) A SECTION HEADING
\vspace*{-0.5pt}
\noindent
In his famous 1949 paper,$^1$ W.~N.~Bailey notes the following
seemingly trivial result.
\begin{lemma}
Let $\alpha=\{\alpha_L\}_{L\geq 0},\ldots,\delta=\{\delta_L\}_{L\geq 0}$
be sequences which satisfy
\begin{equation}
\beta_L = \sum_{k=0}^L \frac{\alpha_k}{(q)_{L-k}(aq)_{L+k}}
\label{ab}
\end{equation}
and
\begin{equation}
\gamma_L = \sum_{k=L}^\infty \frac{\delta_k}{(q)_{k-L}(aq)_{k+L}} \, ,
\label{gd}
\end{equation}
with $(a)_k=(1-a)(1-aq)\cdots (1-aq^{k-1})$.
Then
\begin{equation}
\sum_{L=0}^{\infty} \alpha_L \gamma_L =
\sum_{L=0}^{\infty} \beta_L \delta_L \, .
\label{agbd}
\end{equation}
\end{lemma}
A pair of sequences $(\alpha,\beta)$ which satisfies (\ref{ab})
is said to form a Bailey pair relative to $a$.

Bailey noted that (\ref{agbd}) can be used to obtain identities
of the Rogers--Ramanujan type, provided one finds an appropriate $\delta$
such that (\ref{gd}) can be summed to yield an explicit expression
for $\gamma$.
In particular, taking
\begin{equation}
\delta_L = \frac{q^{L^2} a^L}{(q)_{M-L}} \quad 0 \leq L \leq M 
\qquad {\rm and } \qquad 
\delta_L=0 \quad L > M
\label{delta}
\end{equation}
and applying the $q$-analogue of Saalsch\"utz's theorem, yields
\begin{equation}
\gamma_L = \frac{q^{L^2} a^L}
{(q)_{M-L}(aq)_{M+L}} 
\quad 0 \leq L \leq M 
\qquad {\rm and } \qquad 
\gamma_L=0 \quad L > M.
\label{gamma}
\end{equation}
Substituted into (\ref{agbd}) this gives
\begin{equation}
\sum_{L=0}^M \frac{q^{L^2} a^L}
{(q)_{M-L}(aq)_{M+L}} \; \alpha_L =
\sum_{L=0}^M
\frac{q^{L^2} a^L}{(q)_{M-L}} \; \beta_L \, ,
\label{fin}
\end{equation}
which, after taking $M\to \infty$, simplifies to
\begin{equation}
\frac{1}{(aq)_{\infty}}\sum_{L=0}^{\infty} q^{L^2} a^L \, \alpha_L =
\sum_{L=0}^{\infty}
q^{L^2} a^L \, \beta_L \, .
\label{inf}
\end{equation}

As an example of how this result may be used, we follow
Andrews,$^2$ and take the
following Bailey pair relative to $a$:
\begin{eqnarray}
\alpha_L \!\! &=& \!\! \frac{(-1)^L (1-aq^{2L})(a)_L q^{L(L-1)/2}}
{(1-a)(q)_L} \nonumber \\[2mm]
\beta_L \!\! &=& \!\! \delta_{L,0} \, .
\label{bp}
\end{eqnarray}
Substituting this into (\ref{inf}) and setting $a=q^{\ell}$,
$\ell=0,1,2,\ldots$, we arrive at Euler's identity
\begin{equation}
\frac{1}{(q)_{\infty}} \sum_{j=-\infty}^{\infty} (-1)^j
q^{(3j+1)j/2} = 1,
\end{equation}
independent of $\ell$.

More complicated Bailey pairs were used by Bailey$^1$ and
subsequently by \linebreak Slater$^{3,4}$, who took the Bailey lemma as
starting point for the derivation of her celebrated list of
130 Rogers--Ramanujan identities.

%\pagebreak
%\textheight=7.8truein
%\setcounter{footnote}{0}
%\renewcommand{\thefootnote}{\alph{footnote}}

\section{The Bailey Chain}
\noindent
A particularly important observation was made by Andrews,$^2$ who noted that
if $(\alpha,\beta)$ forms a Bailey pair relative to $a$,
then (\ref{fin}) allows one to construct a new pair of sequences 
$(\alpha',\beta')$ which again forms a Bailey pair relative to $a$.
Specifically, from (\ref{fin}) we infer the following lemma.
\begin{lemma}
If $(\alpha,\beta)$ forms a Bailey pair relative to $a$, then
$(\alpha',\beta')$, defined as 
\begin{eqnarray}
\alpha'_L \!\! &=& \!\!
q^{L^2} a^L \, \alpha_L \nonumber \\[2mm]
\beta'_L \!\! &=& \!\!  \sum_{k=0}^L \frac{q^{k^2} a^k}
{(q)_{L-k} } \; \beta_k ,
\label{Bc}
\end{eqnarray}
again forms a Bailey pair relative to $a$.
\end{lemma}
Since (\ref{Bc}) can of course be iterated an arbitrary number of times,
the above lemma  gives rise to the so-called Bailey chain.$^2$

For example, lemma~2 
applied $k$ times to the Bailey pair (\ref{bp}), gives a new Bailey pair
$(\alpha^{(k)},\beta^{(k)})$ which
substituted into (\ref{inf}) yields
\begin{equation}
\frac{1}{(q)_{\infty}} \sum_{j=-\infty}^{\infty} (-1)^j
q^{\bigl((2k+3)j+1\bigr)j/2} \, a^{jk} =
\sum_{n_1 \geq \ldots \geq n_k \geq 0}
\frac{q^{n_1^2 + \cdots + n_k^2} \, a^{n_1+ \cdots + n_k}}
{(q)_{n_1-n_2} \ldots (q)_{n_{k-1}-n_k} (q)_{n_k}},
\label{AG}
\end{equation}
with $a=1,q$. 
After rewriting the left-hand side into product form
using Jacobi's triple product identity, these identities yield
a subset of Andrews' analytic form of Gordon's identities.$^{5,6}$
For $k=1$ (\ref{AG}) corresponds to the Rogers--Ramanujan 
identities.$^{7,8}$

\section{A Higher-Level Bailey Lemma}
As we have illustrated in the previous sections, Bailey's lemma and
the Bailey chain are extremely powerful concepts, producing an infinite
series of $q$-identities from a single Bailey pair\footnote{By
applying the so-called Bailey lattice,$^{9}$ even larger classes
of identities can be derived from a given Bailey pair.}.
In the following we generalize ($\gamma,\delta$) of (\ref{delta}) and
(\ref{gamma}) to $(\gamma^{(N)},\delta^{(N)})$ where $(\gamma,\delta)
=(\gamma^{(1)},\delta^{(1)})$. This provides, together with the known Bailey
pairs $(\alpha,\beta$) and Bailey's lemma (\ref{agbd}), a vast
number of new $q$-series identities. 
Inserting $(\gamma^{(N)},\delta^{(N)})$ in (\ref{agbd}) gives the 
``higher-level Bailey lemma'', so-called since among many
identities, it gives rise to
Rogers--Ramanujan type identities for the branching functions
of the level-$N$ coset conformal field theories
$({\rm A}^{(1)}_1)_N\times ({\rm A}^{(1)}_1)_L/({\rm A}^{(1)}_1)_{L+N}$.

Before we give our result some more notation is needed.
First, we need the Gaussian polynomial or $q$-binomial coefficient
\begin{equation}
\Mult{A}{B} = \left\{
\begin{array}{ll}
\displaystyle
\frac{(q)_A}{(q)_B(q)_{A-B}}  \quad & {\rm if}  \; 0\leq B \leq A \\[4mm]
0 & \rm otherwise.
\end{array} \right.
\end{equation}
Furthermore, we fix the integer $N\geq 1$, and denote $C$ the
Cartan matrix and ${\cal I}$ the
incidence matrix of the Lie algebra A$_{N-1}$.
That is, ${\cal I}_{j,k} = \delta_{j,k-1}+\delta_{j,k+1}$ and
$C = 2 I - {\cal I}$ with $I$ the $(N-1)\times (N-1)$ identity matrix.
Finally, $\vec{k}$ (where $\vec{k}$ stands for $\vec{m}, \vec{n},
\vec{\mu}, \vec{\eta}$)
and $\vec{\rm e}_{\ell}$ are
$(N-1)$-dimensional vectors with non-negative integer
entries $(\vec{k})_j=k_j$
and $(\vec{\rm e}_{\ell})_j = \delta_{\ell,j}$.

Our result can then be stated as follows.
\begin{lemma}
\label{higherN}
Fix integers $M\geq 0$, $N\geq 1$ and $0\leq \ell \leq N$,
fix $a$ in (\ref{gd}) to $a=q^{\ell}$,
and choose $\delta^{(N)}$ as
\begin{equation}
\delta_L^{(N)} =
\frac{q^{L(L+\ell)/N} }{(q)_{M-L}}
\sum_{\frac{L}{N}-(C^{-1} \vec{n})_1 \in Z}
q^{\vec{n}\, C^{-1} (\vec{n}-\vec{\rm e}_{\ell})}
\prod_{j=1}^{N-1} \Mult{m_j+n_j}{n_j},
\label{deltaN}
\end{equation}
with $0\leq L \leq M$ ($\delta_L=0$ for $L>M$) and
with $\vec{m}$ fixed by $\vec{n}$ through the $(\vec{m},\vec{n})$-system
\begin{equation}
\label{mn}
\vec{m}+\vec{n} = \frac{1}{2}
( {\cal I}\:\vec{m}+
(2L+\ell)\,\vec{\rm e}_{N-1}+\vec{\rm e}_{\ell}).
\end{equation}
Then $\gamma^{(N)}$ is given by
\begin{equation}
\gamma_L^{(N)} = 
\frac{q^{L(L+\ell)/N} }{(q)_{M-L}(q^{\ell+1})_{M+L}}
\sum_{\frac{L}{N}-(C^{-1} \vec{\eta})_1 \in Z}
q^{\vec{\eta}\, C^{-1} (\vec{\eta}-\vec{\rm e}_{\ell})}
\prod_{j=1}^{N-1} \Mult{\mu_j+\eta_j}{\eta_j},
\label{gammaN}
\end{equation}
for $0\leq L \leq M$ ($\gamma_L=0$ for $L>M$),
with $(\vec{\mu},\vec{\eta})$-system
\begin{equation}
\vec{\mu}+\vec{\eta} = \frac{1}{2}
( {\cal I}\:\vec{\mu}+
(M-L)\,\vec{\rm e}_1+(M+L+\ell)\,\vec{\rm e}_{N-1}+\vec{\rm e}_{\ell}).
\label{em}
\end{equation}
\end{lemma}
We note that the sum
$\sum_{A-(C^{-1} \vec{k})_1 \in Z}$
has to be interpreted as a sum over non-negative integers
$k_1,\ldots,k_{N-1}$, such that
$A-\vec{\rm e}_1 C^{-1} \vec{k}$ is again integer.
We further note that for $N=1$ we reproduce the results
(\ref{delta}) and (\ref{gamma}) with $a=1$ or $a=q$.
An inductive proof of lemma \ref{higherN} will be given in Ref.~10.

An immediate corollary of lemma~3 is the generalization of 
(\ref{inf}) to arbitrary $N$.
\begin{corollary}
Let $N$ and $\ell$ be defined as in lemma~3 and let $(\alpha,\beta)$
be a Bailey pair relative to $q^{\ell}$.
Then
\begin{eqnarray}
\label{corol}
\lefteqn{
\frac{1}{(q^{\ell+1})_{\infty}}
\sum_{L=0}^{\infty}
q^{L(L+\ell)/N}  \, \alpha_L
\sum_{\frac{L}{N}-(C^{-1} \vec{\eta})_1 \in Z}
\frac{q^{\vec{\eta}\, C^{-1} (\vec{\eta}-\vec{\rm e}_{\ell})}}
{(q)_{\eta_1} \ldots (q)_{\eta_{N-1}}} } \nonumber \\
&=&
\sum_{L=0}^{\infty}
q^{L(L+\ell)/N} \, \beta_L
\sum_{\frac{L}{N}-(C^{-1} \vec{n})_1 \in Z}
q^{\vec{n}\, C^{-1} (\vec{n}-\vec{\rm e}_{\ell})}
\prod_{j=1}^{N-1} \Mult{m_j+n_j}{n_j},
\label{cor}
\end{eqnarray}
with the $(\vec{m},\vec{n})$-system (\ref{mn}).
\end{corollary}

As a simple application of corollary~1, we substitute the Bailey pair 
(\ref{bp}) with $a=1$ ($\ell=0$) into (\ref{cor}).
This yields the following generalization of Euler's identity:
\begin{equation}
\frac{1}{(q)_{\infty}}
\sum_{j=-\infty}^{\infty}
(-1)^j q^{\bigl((1+2/N)j+1\bigr)j/2}
\sum_{\frac{j}{N}-(C^{-1} \vec{\eta})_1 \in Z}
\frac{q^{\vec{\eta}\, C^{-1} \vec{\eta}}}{
(q)_{\eta_1} \ldots (q)_{\eta_{N-1}}} = 1.
\end{equation}
As a more elaborate example, we substitute the Bailey pair
$(\alpha^{(k)},\beta^{(k)})$ (obtained from (\ref{bp}) by $k$ times 
iterating (\ref{Bc})) into (\ref{cor}).
For $a=1$ ($\ell=0$), this leads to the following generalization of the
first Rogers--Ramanujan ($N=1$, $k=1$) 
and first Andrews--Gordon identity ($N=1$, $k\geq 2$):
$$
\frac{1}{(q)_{\infty}}
\sum_{j=-\infty}^{\infty} (-1)^j
q^{\bigl((2k+1+2/N)j+1\bigr)j/2}
\sum_{\frac{j}{N}-(C^{-1} \vec{\eta})_1 \in Z}
\frac{q^{\vec{\eta}\, C^{-1} \vec{\eta}}}
{(q)_{\eta_1} \ldots (q)_{\eta_{N-1}}}
$$
\vspace{-2mm}
\begin{equation}
= \sum_{r_1 \geq \ldots \geq r_k \geq 0}
\frac{q^{r_1^2/N + r_2^2 + \cdots + r_k^2}} % \, a^{n_1+ \cdots + n_k}}
{(q)_{r_1-r_2} \ldots (q)_{r_{k-1}-r_k} (q)_{r_k}}
\sum_{\frac{r_1}{N}-(C^{-1} \vec{n})_1 \in Z}
q^{\vec{n}\, C^{-1} \vec{n}}
\prod_{j=1}^{N-1} \Mult{m_j+n_j}{n_j},
\label{AGN}
\end{equation}
with $(\vec{m},\vec{n})$-system
\begin{equation}
\vec{m}+\vec{n} = \frac{1}{2}
( {\cal I}\:\vec{m}+
2r_1 \,\vec{\rm e}_{N-1}).
\end{equation}
These identities are closely related to the 
G\"ollnitz--Gordon identities$^{11,12,13}$ when $N=2$.
To the best of our knowledge (\ref{AGN}) is new for $N\geq 3$.

\section{Discussion}
In this note we have presented a higher-level generalization of the 
well-known Bailey lemma. As an application,
some new $q$-series identities of the Rogers--Ramanu\-jan type
have been derived.
%An inductive proof and further applications of our results will
%be given in a future publication.$^10$

Finally, we list some general remarks about the higher-level Bailey
lemma.
\begin{itemize}
\item
Besides the factor $1/(q)_{M-L}(aq)_{M+L}$,
$\gamma_L^{(N)}$ in (\ref{gammaN}) also depends on $M$ through the
$M$-dependence of the $(\vec{\mu},\vec{\eta})$-system (\ref{em}).
Hence, for $N \geq 2$, equations
(\ref{deltaN}) and (\ref{gammaN}) cannot be used to
obtain a Bailey chain in the same way as (\ref{delta}) and (\ref{gamma})
gave (\ref{Bc}). 
\item
Apart from the factor $q^{L(L+\ell)}/(q)_{M-L}$, the expression
(\ref{deltaN})
for $\delta_L^{(N)}$ coincides with the fermionic polynomial
expressions for the configuration sums of
the level-2 A$^{(1)}_{N-1}$ Jimbo--Miwa--Okado models,$^{14}$
as calculated by Foda {\em et al}.$^{15}$
This suggests that other exactly solvable lattice models,
in particular the level-2 ${\cal G}^{(1)}_r$ models with
${\cal G}=\;$D and E,$^{16,17}$ can be used to obtain further generalizations
of the Bailey lemma.
\item
In Ref.~18 it was pointed out that the polynomial identities for
finitized Virasoro characters of the minimal models $M(2,2k+1)^{19,20}$
and $M(p,p+1)^{21-25}$ give rise to Bailey pairs.
The application of Bailey's original lemma to the Bailey pairs arising
from the polynomial identities for the most general model $M(p,p')$ was 
discussed in Ref.~26 and 27 where the Bailey
transformation was interpreted as a renormalization group flow between
different minimal models.
When substituted in (\ref{corol}), the $M(p,p')$ Bailey pairs
yield Rogers--Ramanujan type identities for the branching functions
of the level-$N$ cosets
\begin{equation}
\frac{({\rm A}_1^{(1)})_N \times ({\rm A}_1^{(1)})_{L}}
{({\rm A}_1^{(1)})_{N+L}}\, ,
\label{coset}
\end{equation}
where $N$ is that of lemma~3, and $L$ is
the (in general) fractional level $p'/p-2$ or $-N-2-p'/p$.
For example, from the Bailey pair arising from $M(1,p)$ we arrive at the
unitary character identities with integer level $L=p-2$, previously obtained 
in Refs.~28 and 29. 
A more detailed discussion of the $q$-series identities
for the level-$N$ cosets will be given in Ref.~10.
\item
The original $(\gamma,\delta)$ pair of Bailey$^{1}$ depends on two continuous
parameters $\rho_1$ and $\rho_2$
\begin{equation}
\delta_L = \frac{(\rho_1)_L (\rho_2)_L (aq/\rho_1\rho_2)^L}
{(aq/\rho_1)_M (aq/\rho_2)_M} \:
\frac{(aq/\rho_1\rho_2)_{M-L}}{(q)_{M-L}} \qquad 0 \leq L\leq M
\label{dr}
\end{equation}
and
\begin{equation}
\gamma_L = \frac{(\rho_1)_L (\rho_2)_L (aq/\rho_1\rho_2)^L}
{(aq/\rho_1)_L (aq/\rho_2)_L} \: 
\frac{1}{(q)_{M-L}(aq)_{M+L}} \qquad 0 \leq L\leq M,
\label{gr}
\end{equation}
which reduce to (\ref{delta}) and (\ref{gamma}) as $\rho_1,\rho_2\rightarrow
\infty$. Slater$^{3,4}$ exploits this dependence on $\rho_1$ and $\rho_2$
to obtain many character identities for cosets of the form (\ref{coset})
with $N=2$, by keeping $\rho_1$ finite and letting $\rho_2\rightarrow \infty$.
Further examples of this construction were given in Ref.~27 where it was
also shown that the characters of the unitary $N=2$ supersymmetric models
follow by specializing both $\rho_1$ and $\rho_2$ to appropriate finite
values.
At present it is unclear to us how to generalize the
sequences $(\delta^{(N)},\gamma^{(N)})$ of lemma~3 to include
such additional parameters $\rho_i$.
\item
In Ref.~30 Milne and Lilly have given yet another
generalization of the Bailey lemma by extending the
definitions  (\ref{ab}) and (\ref{gd})
to higher-rank groups ((\ref{ab}) and (\ref{gd}) correspond to A$_1$). 
Using the theory of
higher-rank basic hypergeometric series they then
found an appropriate generalization of $\delta$ in (\ref{dr})
that can be summed explicitly.
An extremely challenging problem would be
to generalize our higher-level
Bailey lemma to the higher-rank cases of Milne and Lilly.
\end{itemize}

\nonumsection{Acknowledgements}
\noindent
It is a pleasure to thank Barry McCoy for many valuable discussions.
SOW thanks the Institute for Theoretical Physics at Stony Brook
for kind hospitality.
This work is (partially) supported by the Australian Research Council 
and NSF grant DMR9404747.

\nonumsection{References}

\end{document}